\documentclass{amsart}

\usepackage{textcomp}
\usepackage{amsmath,amsthm}
\usepackage{amsfonts}
\usepackage{amssymb}
\usepackage{color,hyperref}
\usepackage[table,usenames,dvipsnames]{xcolor}
\usepackage{url}
\usepackage{faktor}
\usepackage{cite}
\usepackage[backgroundcolor=orange!10]{todonotes}

\usepackage{booktabs}
\newcommand{\myrowcolour}{\rowcolor[gray]{0.925}}
\usepackage{txfonts}

\usepackage{tikz}
\usetikzlibrary{matrix,arrows,automata,positioning,arrows.meta,backgrounds,fit,calc}
\tikzset{arr/.style={-Latex}}
\usepackage{pgfmath}
\pgfdeclarelayer{background layer}
\pgfsetlayers{background layer,main}
\usepackage{tikz-cd}

\definecolor{grey}{gray}{0.95}

\definecolor{gray9}{gray}{0.9}
\definecolor{gray8}{gray}{0.8}
\definecolor{gray7}{gray}{0.7}
\definecolor{gray6}{gray}{0.6}

\definecolor{darkblue}{rgb}{0.0,0.0,0.3}
\hypersetup{colorlinks,breaklinks,
            linkcolor=darkblue,urlcolor=darkblue,
            anchorcolor=darkblue,citecolor=darkblue}

\usepackage{orcidlink}

\newcommand{\hooklongrightarrow}{\lhook\joinrel\longrightarrow}

\newcommand{\cT}{{\mathcal T}}
\newcommand{\cS}{{\mathcal S}}

\newcommand{\Z}{\mathbb{Z}}

\newcommand{\BinRel}{\text{B}}
\newcommand{\InvMon}{\mathcal I}
\newcommand{\DualInvMon}{{\mathcal I}^*}
\newcommand{\FullTrans}{\mathcal T}
\newcommand{\Symmetric}{\mathcal S}
\newcommand{\Brauer}{\mathfrak B}
\newcommand{\Partition}{\mathcal P}
\newcommand{\PartBinRel}{\Partition\BinRel}
\newcommand{\TemperleyLieb}{\text{TL}}
\newcommand{\PartialTrans}{\text{P}\FullTrans}

\DeclareMathOperator{\Aut}{Aut}

 \theoremstyle{plain}
 \newtheorem{theorem}{Theorem}[section]

 \theoremstyle{definition}
 \newtheorem{definition}[theorem]{Definition}
 \newtheorem{example}[theorem]{Example}

\newcommand{\SubSemi}{\textsc{SubSemi}}
\newcommand{\GAP}{\textsc{Gap}}

\newcommand{\Clojure}{\textsc{Clojure}}

\newcommand{\Semigroups}{\textsc{Semigroups}}

\newcommand{\grape}{\textsc{GRAPE}}

\begin{document}

\title{Computing  Embeddings and Isomorphisms of Finite Semigroups}
\thanks{This project was funded in part by the Kakenhi grant
22K00015 by the Japan Society for the Promotion of Science (JSPS), titled `On progressing human understanding in the shadow of superhuman
deep learning artificial intelligence entities' (Grant-in-Aid for Scientific
Research type C, \url{https://kaken.nii.ac.jp/grant/KAKENHI-PROJECT-22K00015/}).}

\author{James East}
\address{ Centre for Research in Mathematics and Data Science,
Western Sydney University,
 Sydney, Australia
\textnormal{\orcidlinkf{0000-0001-6112-9754}}}

\author{Attila Egri-Nagy}
\address{Human \& AI Center, Akita International University
 Akita, Japan
\textnormal{\orcidlinkf{0000-0001-7861-7076}}}

\author{Andrew R. Francis}
\address{School of Mathematics and Statistics,
 University of New South Wales,
 Kensington, Australia,
 \textnormal{\orcidlinkf{0000-0002-9938-3499}}}

\author{James D. Mitchell}
\address{Centre for Interdisciplinary Research in Computational Algebra, University of St Andrews, St Andrews, Scotland \textnormal{\orcidlinkf{0000-0002-5489-1617}}}

 \maketitle

\begin{abstract}
Semigroup theory is a branch of abstract algebra, and it provides mathematical tools for the theory of computation.
Finite semigroups can describe state transition systems and thus they model physically realizable computers.
Engineering questions like \emph{What is the minimal number of states to realize a particular computation?}
and
\emph{Which type of computation is more capable?} translate into the algebraic tasks of constructing isomorphisms and embeddings between semigroups of different representations.
The underlying problem is (sub)graph isomorphism, which is computationally difficult in general.
We describe variations of backtrack search algorithms that exploit the algebraic properties of semigroups, and we carry out computational experiments to extend our algebraic knowledge.
In particular, we report new computational results on transformation semigroups and on the more general family of diagram semigroups.
 We study the minimal degree representation problem, count distinct embeddings and work on an open problem of embedding into 2-generated subsemigroups.
\end{abstract}



\section{Introduction}

Semigroups are abstract algebraic structures with a single associative
binary operation, often called multiplication or composition.
Since only the associativity condition is required, semigroup elements
come in a great variety (numbers, functions, matrices, relations,
diagrams, etc.).
Some models of computation, e.g., finite state automata, are based on the composition of total functions and thus automata are in essence transformation semigroups.

Deciding whether or not a semigroup $S$ embeds into another semigroup $T$
($T$ contains a copy of $S$), or whether they are isomorphic ($S$ and $T$ are
essentially the same) is a non-trivial task.
It can be reduced to the graph-isomorphism GI problem of undirected graphs \cite{Mi79}, which is one of the leading problems in theoretical computer science.
It is an open question whether GI can be decided in polynomial time or not \cite{2020_GI_review}.
The more general problem, the subgraph isomorphism problem, is NP-complete \cite{1971CookSubgraphIsomorphism}.

Interpreting semigroups as models of computation, the embedding
problem is the question whether one computer can emulate another one.
This provides a framework for discussing engineering questions like
\emph{What is the minimal number of states needed to realize a
  computation?}, \emph{How can we compare the computational power of
  different forms of computing?} \cite{egri2017finite}.
This practical interest, combined with the inherent difficulty of the problem, justifies an experimental computational approach.

Here we describe a customized
backtrack search algorithm that hugely improves the current
computational capabilities in semigroup theory.
By finding embeddings into full transformation semigroups efficiently, we can solve the minimal degree representation problem up to degree $n=8$ (at least).
We can also find all distinct embeddings.
The computational implementation allowed us to solve open problems about embeddability into 2-generated subsemigroups.

In this section we give the basic definitions of semigroup theory and introduce diagram representations.
Section \ref{sect:backtrack} describes how the backtrack search algorithm can be customized for computing embeddings and deciding isomorphisms. Section \ref{sect:improvements} goes through techniques for improving the basic algorithm. Section \ref{sect:implementations} summarizes existing computational implementations and how the software packages developed here advance the state of the art. Section \ref{sect:results} is about applying the developed algorithms to problems about diagram semigroups.

\subsection{Semigroups -- Basic Definitions and Notation}

Using the mathematical approach, we first give the abstract definition of semigroups and then proceed to the concrete, state-based transformation representation, and finally to the generalized diagram representations.
The interplay between the abstract and the combinatorial representations is central to this paper.

\subsubsection{Semigroups and Homomorphisms}
A \emph{semigroup} is a set $S$ together with an associative binary operation $S\times S\rightarrow S$.
We can denote the binary operation explicitly, like $s\cdot t$, or we can simply use concatenation $st$ for denoting the product of semigroup elements $s,t\in S$.

A \emph{subsemigroup} of $S$ is a subset $U$ that is closed under the operation, denoted by $U\leq S$.
For $A\subseteq S$, $\langle A\rangle$ denotes the least subsemigroup of $S$
containing $A$, the semigroup \emph{generated} by $A$. In other words,
$\langle  A \rangle$ is the semigroup that arises from taking products
of elements from $A$.  If
$\langle  A \rangle=S$, then $A$ is a \emph{generating set} for $S$.

Given two semigroups $(S,\cdot)$ and $(T,\star)$, a
\emph{homomorphism} is a function $\varphi: S\rightarrow T$, such that
$\varphi(u\cdot v)=\varphi(u)\star\varphi(v)$ for all $u,v\in S$.
We call $S$ the \emph{source}, and $T$ the \emph{target} semigroup.
An \emph{isomorphism} is an invertible homomorphism; when an isomorphism exists we write
$S\simeq T$.
Isomorphic structures are essentially the same.
An \emph{embedding} is an injective homomorphism, denoted by
$S\hookrightarrow T$. The image of an embedding of $S$ into $T$ is
an isomorphic copy of $S$ inside $T$.

\subsubsection{Multiplication Tables}
For a finite semigroup $S$ with $|S|=n$, we fix an order of its elements $s_1,\ldots, s_n$, and refer to the elements by their indices.
The  \emph{multiplication table}, or \emph{Cayley table} of $S$
is a matrix $S_{n \times n}$ with entries from
$\{1,\ldots,n\}$, such that $S_{i,j}=k$ if $s_is_j=s_k$, and we can
simply write $ij=k$.
The $i$th row is $S_{i,*}$  and the $j$th column is $S_{*,j}$.
A multiplication table is a complete description of a given semigroup: it contains the result of the multiplication for any pair of elements.
However, it does not contain any internal details of the semigroup elements.
It does not explain why $st=v$, it merely records this fact.
That's why it is an \emph{abstract representation}.

\subsubsection{Transformation Semigroups}
For a finite set $X$, a \emph{transformation} is a function $f:X\rightarrow X$.
If $|X|=n\in\mathbb{N}$, then the set is denoted by positive integers $\{1,\ldots,n\}$ and a transformation $t$ is denoted by simply listing the images of the points: $[t(1),t(2),\ldots,t(n)]$.
A \emph{transformation semigroup} $(X,S)$ of \emph{degree} $n$ is a collection $S$ of transformations of an $n$-element set $X$, closed under function composition.
The semigroup of all  transformations of $n$ points is the  \emph{full transformation semigroup} $\cT_n$.

\subsubsection{Permutation Groups}
\emph{Groups} are semigroups with an identity and unique inverses for all elements.
Groups model reversible, and semigroups model irreversible processes.
The group  consisting of all \emph{permutations} (bijective transformations) of degree $n$ is the
\emph{symmetric group} $\cS_n$.
For permutations we use cycle notation. For instance, $[2,3,1]$ is written as $(1,2,3)$.
The \emph{cyclic} (one-generated) group of
order $n$ is denoted by $\Z_n$.
In computing terms, $\Z_n$ is an $n$-counter.

\subsubsection{Conjugation}
Transformations $s,t\in\cT_n$ are \emph{conjugate} if there exists a
permutation $g\in\cS_n$ such that $t=g^{-1}sg$. In other words, we get $t$
if we relabel the points of $s$ according to a permutation $g$.

Transformations give a concrete representation of semigroups.
From the perspective of automata theory, we can consider the set of states to be $X$ and the semigroup elements the transitions of these states.
In this way, semigroups are abstract representations of automata.

\subsection{Diagram Semigroups}

From the point of view of theoretical computer science and automata theory,
transformation semigroups are the most important representations of semigroups.
However, they are just one special case of \emph{diagram semigroups}, where elements can be composed by adjoining
diagrams.
These are described by fundamental mathematical objects such as relations and functions. 
The original interest came from algebras with a basis whose elements can be multiplied diagrammatically (e.g.~\cite{Brauer1937,Jones1994}) and model various situations in theoretical physics \cite{Martin1994}.

\emph{Partitioned binary relations} are the most general type of diagrams we consider, although historically it was the last to be defined \cite{PartBinRel2013}.
They are directed graphs of binary relations, where the vertices are partitioned into two nonempty sets, the ``interfaces'' for diagram composition. Here we consider partitioning only into equal-sized parts.
For a finite set $X$ a \emph{diagram} is a subset of $(X\cup X')\times (X\cup X')$ where $|X|=|X'|=n$, the \emph{degree} of the diagram, and $X\cap X'=\varnothing$.
Pictorially, we draw the points from $X$ on an upper row with those from $X'$ below, and we draw a directed edge $a\to b$ for each pair $(a,b)$ from the diagram.
The product $\alpha\beta$ of two diagrams $\alpha$ and $\beta$ (on the same set $X$) is calculated as follows.  We first modify $\alpha$ and $\beta$ by changing every lower vertex $x'$ of $\alpha$ and every upper vertex $x$ of $\beta$ to $x''$.  We then stack these modified diagrams together with $\alpha$ above $\beta$ so that the vertices $x''$ are identified in the middle row (there may now be parallel edges in this stacked graph).  Finally, for each $a,b\in X\cup X'$ we include the edge $a\to b$ in $\alpha\beta$ if and only if there is a path from $a$ to $b$ in the stacked graph (as defined above) for which the edges used in the path alternate between the edges of $\alpha$ and the edges of $\beta$.  An example is given in Figure \ref{fig:PartBinRel} (where, for convenience, the edges of $\beta$ are white so that the kinds of paths referred to above are alternating in colour).  This operation is associative, so the set of all diagrams on the set $X$ forms a semigroup.  When $X=\{1,2,\ldots,n\}$, we denote this semigroup by $\PartBinRel_n$.  The identity element of $\PartBinRel_n$ is the diagram containing the edges $x\to x'$ and $x'\to x$ for each $x$. 

\begin{figure}[h]
\begin{center}
\begin{tikzpicture}
\tikzstyle{blackdot}=[draw=black,circle,fill=black,inner sep=1pt]
\tikzstyle{arrow}=[thick,->,>=angle 60]
\node [blackdot] at (0,0) (u1) {};
\node [blackdot,right of=u1] (u2) {};
\node [blackdot,right of=u2] (u3) {};
\node [blackdot,right of=u3] (u4) {};
\node [blackdot,below of=u1] (d1) {};
\node [blackdot,below of=u2] (d2) {};
\node [blackdot,below of=u3] (d3) {};
\node [blackdot,below of=u4] (d4) {};
\draw [arrow] (u1) edge (d3);
\draw [arrow,bend right] (u3) edge (u4);
\draw [arrow] (d1) edge (u2);
\draw [arrow] (d4) edge (u3);
\node at (-.75,-.5) {$\alpha$};
\begin{pgfonlayer}{background layer}
\fill  [gray6] plot (-.3,0) rectangle (3.3,-1);
\end{pgfonlayer}
\begin{scope}[yshift=-1cm]
\tikzstyle{arrow}=[white,thick,->,>=angle 60]
\node [blackdot] at (0,0) (u1) {};
\node [blackdot,right of=u1] (u2) {};
\node [blackdot,right of=u2] (u3) {};
\node [blackdot,right of=u3] (u4) {};
\node [blackdot,below of=u1] (d1) {};
\node [blackdot,below of=u2] (d2) {};
\node [blackdot,below of=u3] (d3) {};
\node [blackdot,below of=u4] (d4) {};
\draw [arrow] (u3) edge (d2);
\draw [arrow,bend right] (u3) edge (u4);
\draw [arrow,bend right] (d2) edge (d1);
\draw [arrow] (d2) edge (u1);
\node at (-.75,-.5) {$\beta$};
\begin{pgfonlayer}{background layer}
\fill  [gray6] plot (-.3,0) rectangle (3.3,-1);
\end{pgfonlayer}
\end{scope}
\begin{scope}[yshift=-3cm]
\tikzstyle{arrow}=[black,thick,->,>=angle 60]
\node [blackdot] at (0,0) (u1) {};
\node [blackdot,right of=u1] (u2) {};
\node [blackdot,right of=u2] (u3) {};
\node [blackdot,right of=u3] (u4) {};
\node [blackdot,below of=u1] (d1) {};
\node [blackdot,below of=u2] (d2) {};
\node [blackdot,below of=u3] (d3) {};
\node [blackdot,below of=u4] (d4) {};
\draw [arrow] (u1) edge (d2);
\draw [arrow,bend right] (u1) edge (u3);
\draw [arrow,bend right] (u3) edge (u4);
\draw [arrow] (d2) edge (u2);
\draw [arrow,bend right] (d2) edge (d1);
\node at (-.75,-.5) {$\alpha\beta$};
\begin{pgfonlayer}{background layer}
\fill  [gray6] plot (-.3,0) rectangle (3.3,-1);
\end{pgfonlayer}
\end{scope}
\end{tikzpicture}
\begin{tikzpicture}
\tikzstyle{plain}=[rounded corners=3pt, draw]
\draw node [plain](Tn) {$\FullTrans_n$};
\draw node [plain,right of=Tn] (In) {$\InvMon_n$};
\draw node [plain,right of=In] (Ins) {$\DualInvMon_n$};
\draw node [plain,right of=Ins] (Br) {$\Brauer_n$};
\draw node [plain,below of=In] (Sn) {$\Symmetric_n$};

\draw node [plain,below right of=Sn] (1n) {$1_n$};
\draw node [plain,below of=Br] (TemperleyLieb) {$TL_n$};
\draw node [plain,above of=Tn] (PTn) {$\PartialTrans_n$};
\draw node [plain,above of=PTn] (Bin) {\BinRel$_n$};
\draw node [plain,above right =5mm of Bin] (PBn) {$\PartBinRel_n$};
\draw node at (3,2) [plain] (PartitionMonoid) {$\Partition_n$};

\draw  (Tn) -- (Sn);
\draw  (In) -- (Sn);
\draw  (Ins) -- (Sn);

\draw  (PartitionMonoid) -- (In);
\draw  (PartitionMonoid) -- (Ins);
\draw  (PartitionMonoid) -- (Tn);
\draw  (PartitionMonoid) -- (Br);

\draw  (TemperleyLieb) -- (1n);
\draw  (Br) -- (TemperleyLieb);
\draw  (Sn) -- (1n);
\draw  (Sn) -- (Br);

\draw (PTn) -- (Tn);
\draw (PTn) -- (In);
\draw (Bin) -- (PTn);
\draw (PBn) -- (Bin);
\draw (PBn) -- (PartitionMonoid);
\end{tikzpicture}
\end{center}
\caption{(left) Composing partitioned binary relations $\alpha$ and $\beta$. The arrows of the product are induced by paths of the stacked diagram with the property that the consecutive arrows have alternating colors. (right) Embeddings of same degree diagram semigroups.}
\label{fig:PartBinRel}
\end{figure}

Other diagram types can be defined by a set of constraints on the set of arrows (direction, number):
$\BinRel_n$ binary relations \cite{plemmons1970},
$\PartialTrans_n$ partial transformation semigroup \cite{ClassicalTransSemigroups2009},
$\Partition_n$ (bi)partition monoid \cite{PartitionAlgebras2005,Martin1994,Jones1994partition},
$\Brauer_n$ Brauer monoid \cite{Brauer1937},
$\Symmetric_n$ symmetric group\cite{CameronPermGroups99, DixonMortimerPermGroups96},
$\FullTrans_n$ full transformation semigroup \cite{ClassicalTransSemigroups2009},
$\InvMon_n$ symmetric inverse monoid \cite{lawson1998inverse},
$\DualInvMon_n$ dual symmetric inverse monoid \cite{DualSymmetricInverse1998},
and $TL_n$  Temperley-Lieb monoid \cite{TemperleyLieb1971}. For the
partial order of diagram semigroups of the same degree see Figure
\ref{fig:PartBinRel}.

\section{Partitioned backtrack}
\label{sect:backtrack}

Classical backtrack search provides a simple algorithm for
constructing semigroup embeddings and isomorphisms, but ignoring algebraic information about the structures involved makes it inefficient.
Therefore, we partition the elements of the target semigroup by their
suitability for being homomorphic images for the elements of the source semigroup.

\subsection{Algorithm for Embedding Multiplication Tables}

Given semigroups $S$ and $T$ such that $|S|\leq|T|$, we would like to
know if we can
construct an embedding $S\hookrightarrow T$, and if so, actually
construct one.
Let $m=|S|$ and $n=|T|$, so we need a map
$\varphi:\{1,\ldots,m\}\rightarrow\{1,\ldots,n\}$ such that
$\varphi$ is injective: $\varphi(i)=\varphi(j)$ implies $i=j$, and it is a homomorphism: 
$\varphi(i)\varphi(j)=\varphi(ij)$, where multiplication on the left hand side is in $T$ and on the right hand side is in $S$.
\begin{example}
For permutation group $\Z_2=\{(), (1,2)\}$ and
transformation semigroup $\cT_2=\{[1,1],[1,2],[2,1],[2,2]\}$,
using lexicographic ordering of the elements, the maps $1\mapsto 2$, $2 \mapsto 3$ give an embedding.
\begin{center}
\begin{tabular}{c|cc}
$\Z_2$ & 1 & 2 \\
\hline
1&1&2\\
2&2&1\\
\end{tabular}
$\hookrightarrow$
\begin{tabular}{c|cccc}
$\cT_2$ & 1 & \textbf{2} &\textbf{3} &4\\
\hline
1&1&1&4&4\\
\textbf{2}&1&\textbf{2}&\textbf{3}&4\\
\textbf{3}&1&\textbf{3}&\textbf{2}&4\\
4&1&4&1&4\\
\end{tabular}
\end{center}
\end{example}
For finding embeddings, the \emph{brute-force} algorithm goes through all possible maps by checking all possible assignments of the
$m$ elements of $S$ to the $n$ elements of $T$. There are
$\frac{n!}{(n-m)!}$ such maps.
By gradually exploiting more information about
$\varphi, S$ and $T$ we can construct increasingly more efficient algorithms.

\subsubsection{Classical Backtrack}
We can observe that any partial $S\to T$ map violating the compatibility condition of homomorphism cannot be fixed by adding more mappings.
A partial non-solution cannot be
extended to a solution, therefore the classical
\emph{backtrack} method \cite{knuth1998} can be  applied.

Let $p$
be a partial solution represented by a sequence of integers, such that
the $i$th element is $\varphi(i)$. So, $p=\left(\varphi(1),\varphi(2),
\ldots, \varphi(l)\right)$ for some $l\leq |S|$.
Being a homomorphism requires that if $ij=k$ (in $S$), then $\varphi(i)\varphi(j)=\varphi(k)$, (in $T$).
However, $k$ and $\varphi(k)$ may not be in the partial solution yet.
If a product $ij=k$ is not in the domain of the partial map yet, i.e.~$k>l$, then the
product $\varphi(i)\varphi(k)$, should also be \emph{undefined}, i.e.~not
being in the sequence $p$.
Now assume that $p$ could be a homomorphism. Then we extend its
sequence by choosing a new $\varphi(l+1)$ from the remaining elements
of $T$ and check whether the homomorphism property is true for products
containing $l+1$. We also have to check for previous undefined
entries, since they may evaluate to $l+1$ and thus become defined.
If the above conditions are true, then we can continue extending $p$
by $l+2$. If not, then according to backtrack, we choose another
$\varphi(l+1)$, or if there is no such candidate, then going back to
$l$ and looking for alternative $\varphi(l)$, and so on.

\subsubsection{Partitioned Backtrack}

We can improve the search algorithm by reducing the number of choices available at each step.
We precompute equivalence relations on $S$ and $T$ such that elements of a
class in $S$ may only be mapped by an embedding to an element of a single class of
$T$.
In other words, we classify elements of $S$ and $T$ by properties
that are invariant under embeddings.
Therefore, when trying to extend a partial solution, we need to look
for candidates only in the corresponding class.
These classes partition the search space, thus we call the method \emph{partitioned backtrack search}.

For abstract semigroups, one such invariant is the semigroup
generalization of the order of the element (isomorphism of  monogenic semigroups).
In a finite semigroup, taking the powers of an element will eventually have a repeated value.

\begin{definition}[Index-period]
  For an element $a$ in a finite semigroup $S$, we have $a^{m+r}=a^m$.
The \emph{index} of $a$ is
the smallest such integer $m\geq 1$ (the exponent of the first repeated power), and the \emph{period}  of $a$ is the smallest such integer $r\geq 1$ (the length of the cycle).
\end{definition}

The index-period pair $(m,r)$ is an invariant under isomorphisms.
For instance, $(1,1)$ means $a^2=a$, i.e., $a$ is an \emph{idempotent}.

\begin{figure}
\begin{center}
  \begin{tikzpicture}[shorten >=1pt, node distance=2cm, on grid, auto,inner sep=2pt]
    \node[draw, circle] (s)   {$s$};
    \node[draw, circle] (s2) [right of=s] {$s^2$};
    \node[draw, circle] (s3) [right of=s2] {$s^3$};
    \node[draw, circle] (s4) [right of=s3] {$s^4$};

    \path[->,every loop/.append style=-{Latex}]
    (s) edge [arr] (s2)
    (s2) edge [arr] (s3)
    (s3) edge [arr,bend left=30] (s4)
    (s4) edge [arr,bend left=30] (s3);

  \end{tikzpicture}
\end{center}
  \caption{Example of a semigroup element with index-period (3,2). The third power $s^3$ is the least power that gets repeated, as $s^5=s^3$. All higher powers alternate between $s^3$ and $s^4$ in a 2-cycle.}
  \label{fig:indexperiod}
\end{figure}

Figure \ref{fig:indexperiod} visualizes the definition of index-period: first there is a linear succession of powers of $s$ that appear once only, then we hit a cycle.
This orbit is the semigroup generated by a single element $S=\langle s\rangle$, called $monogenic$, or $cyclic$ semigroup.
Permutations always have index 1, as bijections are reversible, the first power is bound to be repeated.
However, index 1 does not imply permutations.
As a trivial example, the identity transformation $[1,2,3]$ and the constant map $[1,1,1]$ both have the index-period pair $(1,1)$.

The index-period pairs are invariant under isomorphisms, thus they give a way to classify candidate targets when constructing an isomorphism map.

\begin{definition}[$\cong_\multimap$]
  In a finite semigroup when two elements $s$ and $t$ have the same index-period pairs, we write
  $s\cong_\multimap t$, and $\cong_\multimap$ is an equivalence relation. 
\end{definition}
Partitioning by $\cong_\multimap$ can detect whether an embedding is possible or not:
for each $A\in \faktor{S}{\cong_\multimap}$ the corresponding class $B\in
\faktor{T}{\cong_\multimap}$ (with same index-period as $A$) should satisfy $|A|\leq |B|$. 
It also potentially reduces the search space. Denoting the class in
$\faktor{T}{\cong_\multimap}$ corresponding to the class $A\in
\faktor{S}{\cong_\multimap}$ by $\varphi(A)$, the size of the search space is given by
$$\prod_{A\in \faktor{S}{\cong_\multimap}}\frac{|\varphi(A)|!}{\left(|\varphi(A)|-|A|\right)!},$$
showing that the efficiency depends on the number of the
classes and their cardinalities.

\begin{example}
If $|S|=5$ and $|T|=10$ the search space size is $\frac{10!}{(10-5)!}=30240$. Assuming
that we can partition $S$ into $\cong_\multimap$-classes of size $2,3$ and $T$ into
corresponding classes of size $3,5$ (and a class of size 2 with
index-period not appearing in $S$), then the search space size is
reduced to $\frac{3!}{(3-2)!}\cdot\frac{5!}{(5-3)!}=360$.
\end{example}

After this theoretical example, the natural question is what is the size distribution of the index-period equivalence classes for semigroups.
As an indicative answer, we check a full transformation semigroup. Table \ref{tab:T7ip} shows the size distribution of the equivalence classes of the $7^7=823543$ transformations in $\FullTrans_7$.
The smallest class has $420$, the largest class has $126000$ elements.
Therefore, by itself the index-period is not a strong enough invariant: there are only a few classes, and their sizes are far from uniformly distributed.

\begin{table} 
  \begin{center}
  \begin{tabular}{c|r}
  index-period & \#transformations \\
  \hline
  $( 1 , 12 )$ &  420 \\
 $( 1 , 10 )$ &  504 \\
 $( 1 , 7 )$ &  720 \\
 $( 2 , 6 )$ &  4200 \\
 $( 4 , 3 )$ &  5040 \\
 $( 2 , 5 )$ &  5040 \\
 $( 3 , 4 )$ &  5040 \\
 $( 5 , 2 )$ &  5040 \\
 $( 6 , 1 )$ &  5040 \\
 $( 1 , 1 )$ &  6322 \\
 $( 1 , 5 )$ &  19152 \\
 $( 1 , 6 )$ &  22050 \\
 $( 5 , 1 )$ &  27720 \\
 $( 2 , 4 )$ &  28980 \\
 $( 3 , 3 )$ &  29400 \\
 $( 4 , 2 )$ &  30240 \\
 $( 1 , 4 )$ &  37590 \\
 $( 1 , 3 )$ &  37800 \\
 $( 1 , 2 )$ &  45843 \\
 $( 2 , 1 )$ &  59472 \\
 $( 4 , 1 )$ &  68880 \\
 $( 2 , 3 )$ &  73080 \\
 $( 3 , 2 )$ &  85260 \\
 $( 3 , 1 )$ &  94710 \\
 $( 2 , 2 )$ &  126000
\end{tabular}
\end{center}
\caption{Number of transformations in the index-period equivalence classes of $\FullTrans_7$.}
\label{tab:T7ip}
\end{table}

\subsection{Deciding Isomorphism of Multiplication Tables}

For isomorphism we can define more
invariants and a finer partitioning of elements by using a stronger
equivalence relation.
The invariant properties can be any isomorphism invariant of a multiplication table that do not take into account any information of the
actual ordering of its elements.

A \emph{frequency distribution} takes a multiset and enumerates its distinct elements paired with the number of occurrences of the elements.
For instance, the  frequency distribution of a row vector
$[7,2,1,2,5,5,2,7]$ is $[[1,1],[2,3],[5,2],[7,2]]$, meaning that 1
appears once, 2 appears three times and 5 and 7 twice.
However, we cannot retain the element information as it depends on the sorting of the semigroup elements.
We keep only the sorted frequency values $[1,2,2,3]$.
Sorting is crucial here to decide whether two distributions are the same or not.
For example, the vector $[2,4,2,4]$ has the same frequency
distribution as $[1,3,3,1]$ and $[2,2,4,4]$, which is $[2,2]$ , but $[2,4,4,4]$ has
a different one, namely $[1,3]$.
We can also take the frequency distribution of frequency values, since
it is derived from data containing no information on the ordering of the elements.
We can define invariant properties both on the element and on the
table level. In addition to the index-period values, the \emph{element level
invariants} are the number of occurrences of the element $i$
\begin{enumerate}
\item in the table, \emph{frequency},
\item in the diagonal of the table, \emph{diagonal frequency},
\item in its row $S_{i,*}$, \emph{row frequency},
\item in its column $S_{*,i}$, \emph{column frequency}.
\end{enumerate}
These invariants contain information about the structure of the semigroups (Green's classes).
We put them together in a single aggregated data structure called
the \emph{element profile}, in order to get a finer equivalence relation.
\begin{definition}[$\cong_\boxplus$]
For a finite semigroup $S$, the equivalence relation of having the same element profile is
written as
$ s\cong_\boxplus t$, $s,t\in S$.
\end{definition}

The \emph{table level invariants} can be used to decide the possibility of
an isomorphism. They are frequency distributions of the
\begin{enumerate}
\item elements,
\item diagonal frequencies,
\item column and row frequencies,

\item idempotents (very strong semigroup invariant),
\item element profiles.
\end{enumerate}
\noindent These invariants are
sensitive enough to tell some groups apart, despite their very special nature (Latin square multiplication table, single principal ideal and idempotent).
Diagonal frequencies can tell some groups apart: $\Z_4\mapsto (2,2)$, $\Z_2\times \Z_2\mapsto (4)$, but it assigns $(2,6)$ to both $D_8$ (dihedral group) and $Q_8$ (quaternion group).
In the group case, element profiles
  reduces to the order of elements, but it can distinguish between
  $D_8$ and $Q_8$. However, this invariant fails to detect the
  difference between some direct and semidirect products. For
  instance, $\Z_8\times \Z_2$ and $\Z_8\rtimes \Z_2$ both have 1
  element of order 1, 3 of order 2, 4 of order 4, and 8 of order 8.
\begin{example} Let $S$ be the semigroup generated by transformations
  $[ 2, 1, 1 ], [ 2, 3, 2 ], $ and  $[ 3, 1, 3 ]$. Using only
  index-period values to classify its 15 elements, the size of the
  search space to find isomorphisms to another representation of $S$
  is 2903040. Classifying by element profiles defined above reduces the
  search space size to 768.
\end{example}

\section{Algorithmic improvements}
\label{sect:improvements}

To scale up practical computations we need to apply several
optimization techniques.
These all exploit some particular mathematical properties of the
semigroups.
The incremental nature of homomorphisms allow parallel execution,
symmetries of the semigroups lead to more efficient enumeration, generator sets allow computing homomorphisms `on the fly', and we can speed up computations by implementing custom representations for different diagram semigroups.

\subsection{Parallel Execution}
A search can be executed easily in  parallel if parts of the search space can
be separated in a way that instances of the algorithm never cross the boundaries, often called as `embarrassingly parallel'.
For backtrack, this can be achieved by starting the search from different partial solutions.
Since homomorphisms take subsemigroups to subsemigroups we can attempt to extend
an embedding $\varphi_U: U\hookrightarrow T$ (a partial solution) to
$\varphi_S: S\hookrightarrow T$, if $U$ is a subsemigroup of $S$.
Therefore, we can utilize several processes to calculate all embeddings.
First we find all embeddings $\varphi_U$ sequentially, assuming that
$U$ is small enough to be calculated quickly. Then we can run searches
starting from each different $\varphi_U$ in parallel in order to find all embeddings $\varphi_T$.

Load balancing is a potential issue.
In the extreme case, a single or few partial solutions can take most of the required work, while the other threads finish quickly.
At the scale of the current experiments, this issue was not observed, since checking a single partial solution is quick relative to the large number of them.
Therefore, grouping several of them into a work package mitigates the problem.  

The implementation uses \Clojure's library for reducers based on \textsc{Java}'s Fork-Join framework \cite{JavaFJ}.
The parallel search algorithms are in a separate package called \textsc{Orbit} \cite{orbit}.
The package exposes a single parameter \texttt{*task-size*} to counter potential load balancing issues.

\subsection{Finding distinct embeddings up to conjugation}
It is possible to find an embedding in different ways (different bijections from the source semigroup to one particular image inside the source semigroup), and to find different embeddings (several isomorphic copies of the image in the target semigroup).
The first case is when the source semigroup has symmetries, i.e.~its
automorphism group $\Aut(S)=\left\{\sigma:S\rightarrow
  S\mid\sigma(S)\simeq S\right\}$ is non-trivial.
Similarly, some pairs of embeddings can be transformed into each other by those symmetry operations.
Both correspond to relabelings of the elements.
We say that those solutions are \emph{conjugate} to each other and
being conjugate is an equivalence relation. Therefore, it is enough to produce a
single representative element from each conjugacy class.

The representative of a conjugacy class of embeddings can be chosen to
be the minimal one in lexicographic ordering.
For a single transformation (endofunction) we can compute the minimal representative in $O(n^2)$ time, where $n$ is the degree of the transformation \cite{2024minrependofunc}.
However, we also have to calculate conjugacy class representatives for sequences and sets of transformations (and other diagram elements).
Those involve enumerating permutations.
The worst case is the complete enumeration of the symmetric group $\Symmetric_n$, which happens for instance for the set of $n$ constant maps of degree $n$.

For sequences, since a partial solution is a sequence, we can easily home in on the
minimal one (as opposed to finding minimal elements for set-wise conjugation).
As a sequence grows, the set of  symmetries which can take the
sequence to its representative is shrinking, since there are more and
more constraints on how the sequence can be transformed into the
corresponding minimal one.
Computational experiments show that the symmetry classes of sequences become singletons long before the embedding is fully defined.
Once we have only one possible such symmetry, we can be sure that any
extension to a full solution will be a conjugacy class representative.

\subsection{Generator tables}

When computing with algebraic structures, as much it is possible, we want to work with the generators only, not the complete structure.
An extreme example is $\FullTrans_n$, the full transformation semigroup of degree $n$, which has $n^n$ elements, but the standard generating set has only 3 elements: a cycle, a transposition and a collapser (see Figure \ref{fig:T3T4embeds} for these generator set).
The fully computed multiplication table would have $n^{2n}$ entries.

For constructing morphisms, it is enough to have the columns corresponding to multiplication by the generators.
In the source semigroup, we can restrict to multiply only by the generators and we map only these generators.
This trick does not work in the target semigroup, as we keep choosing new candidate generators.
However, it is enough to enumerate the suitable
elements only (e.g., the matching index-period equivalence classes).
This works well for `full' structures of certain types where
we can combinatorially enumerate semigroup elements with the required properties instead of
generating  all elements and filtering the candidates.

In a sense, we are building an isomorphism of the Cayley-graphs of the source
and of the image subsemigroup in the target semigroup.
Once we choose a candidate target generator, we can build both semigroups at the same time to see whether multiplications are compatible.
A partial solution for this method is not a sequence (instead of an
array we have a hashtable), thus to enumerate the morphisms up to
relabelings, we need to find minimal class representatives by set-wise
conjugation.

\subsection{Custom Representations}

\begin{figure}
\begin{center}
  \begin{tikzpicture}
\tikzstyle{blackdot}=[draw=black,circle,fill=black,inner sep=1pt]
\tikzstyle{arrow}=[thick]
\node [blackdot] at (0,0) (u1) {};
\node [blackdot,right of=u1] (u2) {};
\node [blackdot,right of=u2] (u3) {};
\node [blackdot,right of=u3] (u4) {};
\node [blackdot,below of=u1] (d1) {};
\node [blackdot,below of=u2] (d2) {};
\node [blackdot,below of=u3] (d3) {};
\node [blackdot,below of=u4] (d4) {};
\draw [arrow] (u1) edge (d2);
\draw [arrow] (d1) edge (u2);
\draw [arrow,bend right] (u3) edge (u4);
\draw [arrow,bend right] (d4) edge (d3);
\node at (-.75,-.5) {$\alpha$};
\begin{pgfonlayer}{background layer}
\fill  [gray6] plot (-.3,0) rectangle (3.3,-1);
\end{pgfonlayer}
\begin{scope}[yshift=-1cm]
\tikzstyle{arrow}=[white,thick]
\node [blackdot] at (0,0) (u1) {};
\node [blackdot,right of=u1] (u2) {};
\node [blackdot,right of=u2] (u3) {};
\node [blackdot,right of=u3] (u4) {};
\node [blackdot,below of=u1] (d1) {};
\node [blackdot,below of=u2] (d2) {};
\node [blackdot,below of=u3] (d3) {};
\node [blackdot,below of=u4] (d4) {};
\draw [arrow,bend right] (u3) edge (u4);
\draw [arrow,bend right] (d2) edge (d1);
\draw [arrow] (d3) edge (u1);
\draw [arrow] (d4) edge (u2);
\node at (-.75,-.5) {$\beta$};
\begin{pgfonlayer}{background layer}
\fill  [gray6] plot (-.3,0) rectangle (3.3,-1);
\end{pgfonlayer}
\end{scope}
\begin{scope}[yshift=-3cm]
\tikzstyle{arrow}=[black,thick]
\node [blackdot] at (0,0) (u1) {};
\node [blackdot,right of=u1] (u2) {};
\node [blackdot,right of=u2] (u3) {};
\node [blackdot,right of=u3] (u4) {};
\node [blackdot,below of=u1] (d1) {};
\node [blackdot,below of=u2] (d2) {};
\node [blackdot,below of=u3] (d3) {};
\node [blackdot,below of=u4] (d4) {};
\draw [arrow,bend right] (u3) edge (u4);
\draw [arrow] (u1) edge (d4);
\draw [arrow] (u2) edge (d3);
\draw [arrow,bend right] (d2) edge (d1);
\node at (-.75,-.5) {$\alpha\beta$};
\begin{pgfonlayer}{background layer}
\fill  [gray6] plot (-.3,0) rectangle (3.3,-1);
\end{pgfonlayer}
\end{scope}
\end{tikzpicture}
\end{center}
\caption{Multiplication in the Brauer monoid $\Brauer_4$. In our custom representation $\alpha=[6,5,4,3,2,1,8,7]$, $\beta=[7,8,4,3,6,5,1,2]$ and the product is $\alpha\beta=[8,7,4,3,6,5,2,1]$. `Loops' formed in the middle are ignored. In general, we may need to trace long paths of alternating black and white edges.}
\label{fig:BrauerMult}
\end{figure}

When computing with semigroups, multiplication is the repeated operation.
Therefore, the embedding algorithms can be improved by making the multiplication faster.
For instance, for diagram semigroups, we can implement a single multiplication algorithm, the most abstract for partitioned binary relations.
In those, for each point we have a \emph{set} of images.
Obviously, this would be a wasteful representation for transformations and permutations, where each point has only a single image, and computing with sets is an unnecessary overhead.
Therefore, computer algebra software packages have custom representations for transformations.

A less obvious case is when we have some other diagram representation, closer to the most general.
In this project we implemented a specific data structure and multiplication for Brauer monoid elements.
The elements are diagrams with $2n$ points, $n$ at the `top', $n$ at the `bottom'.
A point can be connected to exactly one other point.
We represent such a diagram $b$ with a list of $2n$ entries, $[b(1),\ldots,b(2n)]$, with the constraint that $b(i)=j\implies b(j)=i$.
There is no need to enforce these constraints when working with a generating set, since multiplication of two such diagrams produce a third with the same properties.
See Figure \ref{fig:BrauerMult} for an example.

\section{Software Implementations}
\label{sect:implementations}

In computational semigroup theory, the \Semigroups~package \cite{Semigroups} for the \GAP~computer algebra system \cite{GAP4} emerged as the leading solution.
It also serves a base for more specialized packages.

The \SubSemi~package \cite{subsemi} contains the implementations of the above algorithms.
Despite being inefficient for large semigroups due to
the multiplication table representation, they proved to be immensely useful
in practical computations.
The original goal was to extended results of transformation semigroup enumeration \cite{T3enum1970,T3enum1991}.
\SubSemi~enumerated all 4-state finite computations up to isomorphism
\cite{T4enum}.
Moreover, the enumeration was  extended for more general diagram
semigroups \cite{diagsgps}.
Deciding isomorphism was also needed for enumerating independent generating sets of symmetric groups \cite{SnIS}.
These applications required to develop the isomorphism and embedding construction algorithms described in this paper.

For the isomorphism calculation, we can compare our method to two other methods in the 
\Semigroups~package~\cite{Semigroups}.
The older method uses the \texttt{SmallestMultiplicationTable} function  depending on the
\texttt{SmallestImageSet} function in \GAP~\cite{GAP4} and the \grape~package \cite{GRAPEpaper}.
This function calculates the smallest multiplication table in
lexicographic ordering for the semigroup (by calculating a group action orbit of the
symmetric group), which can be used for
isomorphism testing.
Its complexity depends on the size of the semigroup since we have to
act on the table by the corresponding symmetric group, while our method
does not require group actions and it uses more information about
the semigroup structure, thus it is more scalable. In practice, roughly speaking, with
\SubSemi~we can calculate embeddings and isomorphisms as long as the fully enumerated semigroup fits into the memory of a single computer
(e.g.~embeddings into $\FullTrans_8$ with $8^8=16777216$ elements).

The newer method uses graph isomorphism solvers directly for computing isomorphism: the multiplication tables are turned into  vertex coloured (non-directed) graphs, and then checked whether they are isomorphic using \textsc{bliss}\cite{bliss} or \textsc{nauty}\cite{nauty} (depending on which is enabled in the \textsc{Digraphs} package \cite{Digraphs}).
The graph  constructed has many more nodes than the semigroup has elements (they are used to indicate the direction of the edges, and the element that we are multiplying by).

For computational experiments, where the results are not given by mathematical proofs, \emph{verification} by \emph{recomputation} is important.
We routinely reproduced results by other implementations, and developed a `redundant' package for the algorithms presented here.

We recomputed the results of
\cite{ComputingAutGrpOfSgps2010}, using our method on small semigroups \cite{smallsemi}.
We also used this technique for finding isomorphism classes and
determining the automorphism groups of all transformation semigroups
of degree 4 \cite{T4enum}.

For further verification a `shadow' implementation was also developed, \textsc{Kigen} \cite{kigen}.
It is built on a different architecture: \Clojure \cite{2020Clojure}, a purely functional programming language on top of \textsc{Java}, while $\GAP$ is a functional object-oriented language built around a \textsc{C}/\textsc{C++} kernel.
The \textsc{Kigen} package benefits from the seamless parallelization features of \Clojure, a crucial ingredient for the results below.

\section{Computational Experiments and Results}
\label{sect:results}

\subsubsection{Minimal degree realization of computations}
A basic result of semigroup theory is that any semigroup with  $n$ elements can be realized as a
transformation semigroup of at most degree $n+1$ (the semigroup analogue of
Cayley's Theorem for groups, e.g.~\cite{Howie95}, Chapter 1).
In practice, we would like to have a smaller degree transformation
representation.
This is a more general problem than finite automata minimization
(where to goal is to recognize a language efficiently).
This problem naturally generalizes to all types of diagram semigroups.
Given an abstract semigroup, we would like to find minimal diagram
semigroup realizations of it (in terms of the points in the diagrams)
for a given type of diagram. The general strategy is trying to find
embeddings into the degree $n$ full diagram semigroups (containing all
degree $n$ diagrams of that type).
\begin{definition} For a semigroup $S$ the minimal $D$-diagram
  representation degree is
$\mu_D(S)=\min\{n\mid S\hookrightarrow D_n\}$ where $D\in\{\PartBinRel, \BinRel, \PartialTrans, \FullTrans, \Symmetric, \InvMon, \DualInvMon, \Partition, \Brauer, \TemperleyLieb\}$.
\end{definition}

\begin{figure*}
  $$S=\left\{\begin{pmatrix}0&0\\0&0\end{pmatrix},
\begin{pmatrix}0&0\\0&1\end{pmatrix},
\begin{pmatrix}0&0\\1&0\end{pmatrix},
\begin{pmatrix}0&1\\0&1\end{pmatrix},
\begin{pmatrix}1&0\\1&0\end{pmatrix}\right\},\ \ \begin{tabular}{c|ccccc}
$S$ & 1 &2 &3 &4 &5\\
\hline
1&  1&  1&  1&  1&  1 \\
2&    1&  2&  3&  2&  3 \\
3&    1&  1&  1&  2&  3 \\
4&    1&  4&  5&  4&  5 \\
5&    1&  1&  1&  4&  5
\end{tabular}
$$
\caption{An example matrix semigroup with listed elements and the corresponding abstract multiplication table. The numbers in the table correspond to the order of the listed matrices.}
\label{fig:matrix}
\end{figure*}

\begin{example} Finding minimal degree diagram representations of
  a matrix semigroup over $\Z$ by using its multiplication table (Fig.~\ref{fig:matrix}).
The transformation and partition representation both require 3 points $\mu_\FullTrans(S)=\mu_\Partition(S)=3$.
We need two more points to represent $S$ as Brauer and Temperley-Lieb monoid: $\mu_\Brauer(S)=\mu_\TemperleyLieb(S)=5$.
\end{example}

\subsubsection{Comparing `computational power'
} Various diagram representations can be considered as different models of
computation. Traditionally, we compute by composing functions, but we can use binary relations or equivalence relations.
The interesting question is how much `bigger' in size
 a less capable computational device should be in order to realize a
 more powerful one. To what extent can size make up for structure?

For diagram semigroups we can use the minimal degree diagram
representation search to assess relative computational power.
 For instance, to emulate computations of
 equivalence classes with transformations, or $n=2$
 we need at least 7 states: $\Partition_2\hookrightarrow \FullTrans_7$.
 In this case we know precisely the minimum degrees \cite{2025transformationrepresentationsdiagrammonoids}.

The degree 1 partition binary monoid has only 16 elements, but
realized as transformations it requires 6 points:
$\PartBinRel_1\hookrightarrow\FullTrans_6$. This value is
substantially lower
then 17 required by the right regular representation.

Thinking backwards, transformations realized by some stronger
computational model, we find that the distinction between reversible
and irreversible computation gives a hard limit. Despite its huge size
$|\PartBinRel_2|=65536$, partitioned binary relation monoids of lower
degree cannot emulate full transformations semigroups, due to the lack
of permutations.

Here is a list of similar results:
$\Brauer_1\cong \FullTrans_1$, $\Brauer_2\hookrightarrow\FullTrans_3$,
$\TemperleyLieb_1\cong \FullTrans_1$,
$\TemperleyLieb_2\hookrightarrow\FullTrans_2$,
$\TemperleyLieb_3\hookrightarrow\FullTrans_4$, 
$\Partition_1\hookrightarrow\Brauer_2$,
$\Brauer_3\hookrightarrow\BinRel_3$,
$\DualInvMon_2\hookrightarrow\BinRel_3$.
Many of these involve quite heavy computations.

An interesting case is the embedding of the full transformation semigroup into the Brauer monoid
$\FullTrans_2\hookrightarrow\Brauer_3$ but $\FullTrans_3$ does not
embed into $\Brauer_7$.
\label{lab:Tn_Br}
This raised the question whether it would be possible to embed it into $\Brauer_8$.
At this point we needed parallelization to check that it does not.
The embedding already fails for the semigroup generated by the collapser and the cycle of the standard generating sets (see Figure \ref{fig:T3T4embeds} for the generator types).
Missing the transposition generator yields a subsemigroup of $\FullTrans_3$ with 24 elements, allowing a bit faster embedding checks.

\subsubsection{Counting distinct embeddings}
\begin{figure}[ht]
  \begin{center}
  \includegraphics{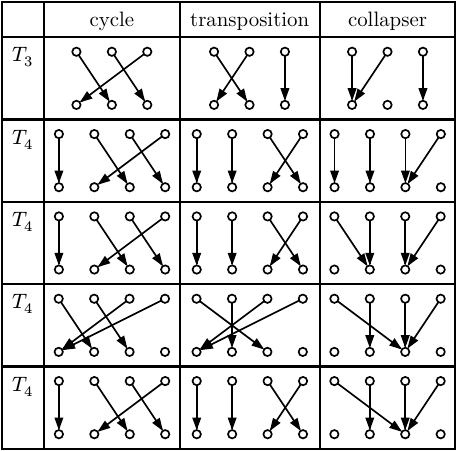}
  \end{center}
  \caption{All 4 (up to conjugation) embeddings of $\FullTrans_3$ into $\FullTrans_4$  given by the mappings of the standard generators. The embeddings simply fix the additional point, or collapse it into
another point. For higher degree embeddings, avoiding the interference with the multiplication in $\FullTrans_3$ yields combinatorial structures of similar fixing and collapsing.}
\label{fig:T3T4embeds}
\end{figure}

Beyond the question whether embedding is possible, we are also interested in
the number of different embeddings, up to conjugation and up to automorphism.
We use the notation $S\overset{k}{\hooklongrightarrow} T$ to indicate that there
are $k$ distinct copies of $S$ inside $T$.
We systematically explored the embeddings of symmetric groups into larger symmetric groups (Table \ref{tab:Snembeddings}).
Similarly, Table \ref{tab:Tnembeddings} summarizes $k$-embeddings for low degree
transformation semigroups.
Figure \ref{fig:T3T4embeds} indicates how these embeddings work.

\begin{table*}[h]
  \begin{center}
    \begin{tabular}{crrrrrrrrr}
      \toprule
      $\Symmetric_m\overset{k}{\hooklongrightarrow} \Symmetric_n$ & $n=1$ & $n=2$& $n=3$& $n=4$& $n=5$&$n=6$&$n=7$ & $n=8$ & $n=9$\\
      \midrule
      $m=1$ & 1 &  1 & 1 & 1 & 1 & 1 & 1 & 1 & 1\\
      \myrowcolour $m=2$&  & 1 &  1 & 2 & 2 & 3 & 3 & 4 & 4\\
      $m=3$ &  &   &  1 & 1 & 2 & 4 & 5 & 7 & 10\\
      \myrowcolour $m=4$ &  &  &   & 1 & 1 & 4 & 5 & 10 & 13\\
      $m=5$ &  &  &   &  & 1 & 2 & 3 & 4 & 5\\
      \myrowcolour $m=6$ &  &  &   &  &  & 1 & 1 & 2 & 2\\
      $m=7$ &  &  &   &  &  &  & 1 & 1 & 2 \\
      \myrowcolour $m=8$ &  &  &   &  &  &  &  & 1 & 1\\
      $m=9$ &  &  &   &  &  &  &  & &  1 \\
      \bottomrule
    \end{tabular}
  \end{center}
  \caption{Number of distinct embeddings of symmetric groups.}
  \label{tab:Snembeddings}
\end{table*}

\begin{table*}
    \begin{center}
    \begin{tabular}{crrrrrrrr}
      \toprule
      $\FullTrans_m\overset{k}{\hooklongrightarrow} \FullTrans_n$ & $n=1$ & $n=2$& $n=3$& $n=4$& $n=5$&$n=6$&$n=7$ & $n=8$ \\
      \midrule
      $m=1$ & 1 &  2 & 3 & 5 & 7 & 11 & 15 & 22\\
      \myrowcolour $m=2$&  & 1 &  3 & 12 & 35 & 110& 309 & 879\\
      $m=3$ &  &   &  1 & 4 & 17 & 64& 221 & 736\\
      \myrowcolour $m=4$ &  &  &   & 1 & 2 & 6 & 16 & 48\\
      $m=5$ &  &  &   &  & 1 & 2 & 6 & 16\\
      \myrowcolour $m=6$ &  &  &   &  &  & 1 & 2 & 6\\
      $m=7$ &  &  &   &  &  &  & 1 & 2\\
      \myrowcolour $m=8$ &  &  &   &  &  &  &  & 1 \\
      \bottomrule
    \end{tabular}
 \end{center}
  \caption{Number of distinct embeddings of full transformation
    semigroups. Embedding the trivial monoid is equivalent to finding
    idempotent elements ($e^2=e$) up to conjugation in the target
    semigroup. An interesting pattern seems to emerge for the
    $\FullTrans_i\overset{k}{\hooklongrightarrow} \FullTrans_n$,
    $i=\{n-3,n-2,n-1\}$ cases.}
  \label{tab:Tnembeddings}

  \end{table*}

There are 282 non-empty transformation semigroups on 3 points up to
conjugation. For degree $4$, this number is  132069775 \cite{T4enum}.
How many isomorphic copies of the degree 3 transformation semigroups
can we find inside $\cT_4$? Counting the embeddings up to conjugation,
we find only 2347 subsemigroups of $\cT_4$ that are isomorphic to some subsemigroup
of $\cT_3$; so most degree 4 transformation semigroups are
`new'. Calculating the same number for $\cT_5$ yields 18236; a modest
increase compared to the still unknown, but expected-to-be gigantic,
number of degree 5 transformation semigroups.

\subsubsection{Embeddings into 2-generated subsemigroups}

Here we will answer a couple of open questions stated in \cite{2014JEPartitionMonoidEmbeddings} using the embedding algorithm.
There the primary interest is in the embeddability into a 2-generated semigroup, denoted by $S\overset{2\text{-gen}}{\hooklongrightarrow}T$.
Being $n$-generated means that the semigroup can be generated by $n$
elements but no less.
It is easy to prove that
$\Brauer_n\overset{2\text{-gen}}{\hooklongrightarrow}\Brauer_{n+2}$,
but there is no 2-generated subsemigroup of $\Brauer_{n+1}$ where we
can embed $\Brauer_n$.
The same was conjectured for the partition monoid, but contrary to the
expectation we found
that
$\Partition_2\overset{2\text{-gen}}{\hooklongrightarrow} \Partition_3$, in
3 different ways.

To obtain this result, we enumerated conjugacy class representatives of subsemigroups of the
target semigroup that are
2-generated, then filtered them for the property of being at least as
big as the source semigroup.
Finally, we used our embedding search to check each of the candidate
2-generated subsemigroups.

\section{Conclusion and Future Work}

We carried out computational experiments to extend our knowledge about diagram semigroups, which are different combinatorial representations of semigroups.
For constructing embeddings and deciding isomorphisms the partitioned
backtrack algorithm reduces the search space by exploiting information
about the source and the target semigroups. Further algorithmic
optimizations improved scalability, thus we could solve open problems
by computational means.

We can extend the optimizations by considering the order of the semigroup elements, borrowing ideas from graph colouring.
The so-called ``Welsh-Powell'' ordering of the nodes (from highest to lowest degree) \cite{WelshPowell67} uses the idea that the definitions made at the start of the sequence imply more restrictions on the later values if they are more connected.
The opposite example shows why this could be useful.
If there is an isolated vertex, and we define its $\varphi$ value first, then this implies no restrictions at all on any subsequent definitions.

Another future task is to do a comprehensive comparison of the now available different methods, especially the performance differences between the specialized and the general search methods.

\subsection*{Software Tools}
The scripts required to reproduce these results are bundled with the software packages \cite{subsemi,kigen}.
The subfolder \texttt{2025\_Computing\_\-Embeddings\_\-and\_\-Isomorphisms\_of} \texttt{\_Finite\_Semigroups} containing the computations for this paper can be found in the \texttt{experiments} folder of the package in both projects.

\bibliographystyle{IEEEtran}
\bibliography{../compsemi.bib}

\begin{thebibliography}{10}
\providecommand{\url}[1]{#1}
\csname url@samestyle\endcsname
\providecommand{\newblock}{\relax}
\providecommand{\bibinfo}[2]{#2}
\providecommand{\BIBentrySTDinterwordspacing}{\spaceskip=0pt\relax}
\providecommand{\BIBentryALTinterwordstretchfactor}{4}
\providecommand{\BIBentryALTinterwordspacing}{\spaceskip=\fontdimen2\font plus
\BIBentryALTinterwordstretchfactor\fontdimen3\font minus
  \fontdimen4\font\relax}
\providecommand{\BIBforeignlanguage}[2]{{%
\expandafter\ifx\csname l@#1\endcsname\relax
\typeout{** WARNING: IEEEtran.bst: No hyphenation pattern has been}%
\typeout{** loaded for the language `#1'. Using the pattern for}%
\typeout{** the default language instead.}%
\else
\language=\csname l@#1\endcsname
\fi
#2}}
\providecommand{\BIBdecl}{\relax}
\BIBdecl

\bibitem{Mi79}
G.~L. Miller, ``Graph isomorphism, general remarks,'' \emph{Journal of Computer
  and System Sciences}, vol.~18, no.~2, pp. 128--142, April 1979.

\bibitem{2020_GI_review}
\BIBentryALTinterwordspacing
M.~Grohe and P.~Schweitzer, ``The graph isomorphism problem,'' \emph{Commun.
  ACM}, vol.~63, no.~11, p. 128–134, 2020. [Online]. Available:
  \url{https://doi.org/10.1145/3372123}
\BIBentrySTDinterwordspacing

\bibitem{1971CookSubgraphIsomorphism}
\BIBentryALTinterwordspacing
S.~A. Cook, ``The complexity of theorem-proving procedures,'' in
  \emph{Proceedings of the Third Annual ACM Symposium on Theory of Computing},
  ser. STOC '71.\hskip 1em plus 0.5em minus 0.4em\relax New York, NY, USA:
  Association for Computing Machinery, 1971, pp. 151--158. [Online]. Available:
  \url{https://doi.org/10.1145/800157.805047}
\BIBentrySTDinterwordspacing

\bibitem{egri2017finite}
A.~Egri-Nagy, ``Finite computational structures and implementations: Semigroups
  and morphic relations,'' \emph{International Journal of Networking and
  Computing}, vol.~7, no.~2, pp. 318--335, 2017.

\bibitem{Brauer1937}
R.~Brauer, ``On algebras which are connected with the semisimple continuous
  groups,'' \emph{Annals of Mathematics}, vol.~38, no.~4, pp. 857--872, 1937.

\bibitem{Jones1994}
V.~F.~R. Jones, ``A quotient of the affine {H}ecke algebra in the {B}rauer
  algebra,'' \emph{Enseign. Math. (2)}, vol.~40, no. 3-4, pp. 313--344, 1994.

\bibitem{Martin1994}
\BIBentryALTinterwordspacing
P.~Martin, ``Temperley-{L}ieb algebras for nonplanar statistical
  mechanics---the partition algebra construction,'' \emph{J. Knot Theory
  Ramifications}, vol.~3, no.~1, pp. 51--82, 1994. [Online]. Available:
  \url{http://dx.doi.org/10.1142/S0218216594000071}
\BIBentrySTDinterwordspacing

\bibitem{PartBinRel2013}
P.~Martin and V.~Mazorchuk, ``Partitioned binary relations,'' \emph{Mathematica
  Scandinavica}, vol. 113, no.~1, pp. 30--52, 2013.

\bibitem{plemmons1970}
R.~J. Plemmons and M.~T. West, ``On the semigroup of binary relations.''
  \emph{Pacific Journal of Mathematics}, vol.~35, no.~3, pp. 743--753, 1970.

\bibitem{ClassicalTransSemigroups2009}
O.~Ganyushkin and V.~Mazorchuk, \emph{{Classical Transformation Semigroups}},
  ser. {Algebra and Applications}.\hskip 1em plus 0.5em minus 0.4em\relax
  Springer, 2009.

\bibitem{PartitionAlgebras2005}
T.~Halverson and A.~Ram, ``Partition {A}lgebras,'' \emph{European J. Combin.},
  vol.~26, no.~6, pp. 869--921, 2005.

\bibitem{Jones1994partition}
V.~F.~R. Jones, ``The {P}otts model and the symmetric group,'' in
  \emph{Subfactors ({K}yuzeso, 1993)}.\hskip 1em plus 0.5em minus 0.4em\relax
  World Sci. Publ., River Edge, NJ, 1994, pp. 259--267.

\bibitem{CameronPermGroups99}
P.~J. Cameron, \emph{{Permutation Groups}}.\hskip 1em plus 0.5em minus
  0.4em\relax London Mathematical Society, 1999.

\bibitem{DixonMortimerPermGroups96}
J.~D. Dixon and B.~Mortimer, \emph{{Permutation Groups}}, ser. {Graduate Texts
  in Mathematics 163}.\hskip 1em plus 0.5em minus 0.4em\relax Springer, 1996.

\bibitem{lawson1998inverse}
M.~Lawson, \emph{Inverse Semigroups: The Theory of Partial Symmetries}.\hskip
  1em plus 0.5em minus 0.4em\relax World Scientific, 1998.

\bibitem{DualSymmetricInverse1998}
\BIBentryALTinterwordspacing
D.~G. Fitzgerald and J.~Leech, ``Dual symmetric inverse monoids and
  representation theory,'' \emph{Journal of the Australian Mathematical Society
  (Series A)}, vol.~64, pp. 345--367, 6 1998. [Online]. Available:
  \url{http://journals.cambridge.org/article_S1446788700039227}
\BIBentrySTDinterwordspacing

\bibitem{TemperleyLieb1971}
E.~H.~L. H.~N. V.~Temperley, ``Relations between the 'percolation' and
  'colouring' problem and other graph-theoretical problems associated with
  regular planar lattices: Some exact results for the 'percolation' problem,''
  \emph{Proc. Roy. Soc. A, Mathematical and Physical Sciences}, vol. 322, no.
  1549, pp. 251--280, 1971.

\bibitem{knuth1998}
D.~Knuth, \emph{{The Art of Computer Programming}}.\hskip 1em plus 0.5em minus
  0.4em\relax Boston, MA, USA: Addison-Wesley Longman Publishing Co., Inc.,
  1998, vol. 1-3.

\bibitem{JavaFJ}
\BIBentryALTinterwordspacing
D.~Lea, ``A {J}ava fork/join framework,'' in \emph{Proceedings of the ACM 2000
  Conference on Java Grande}, ser. JAVA '00.\hskip 1em plus 0.5em minus
  0.4em\relax New York, NY, USA: ACM, 2000, pp. 36--43. [Online]. Available:
  \url{http://doi.acm.org/10.1145/337449.337465}
\BIBentrySTDinterwordspacing

\bibitem{orbit}
A.~Egri-Nagy, ``{\textsc{orbit} Generic sequential and parallel search
  algorithms},''
  \href{https://codeberg.org/egri-nagy/orbit}{\texttt{codeberg.org/egri-nagy/orbit}},
  2025.

\bibitem{2024minrependofunc}
\BIBentryALTinterwordspacing
J.~D. Mitchell, S.~Mukherjee, and P.~Vojt{\v e}chovsk{\'y}, ``Minimal
  representatives of endofunctions,'' \emph{Semigroup Forum}, vol. 109, no.~3,
  pp. 626--638, 2024. [Online]. Available:
  \url{https://doi.org/10.1007/s00233-024-10472-4}
\BIBentrySTDinterwordspacing

\bibitem{Semigroups}
\BIBentryALTinterwordspacing
J.~D. Mitchell \emph{et~al.}, \emph{Semigroups - GAP package, Version 5.5.1},
  Jun 2025. [Online]. Available: \url{http://dx.doi.org/10.5281/zenodo.592893}
\BIBentrySTDinterwordspacing

\bibitem{GAP4}
\emph{{GAP -- Groups, Algorithms, and Programming, Version 4.15.1}}, The
  GAP~Group, 2025, \href{https://www.gap-system.org}{\url{gap-system.org}}.

\bibitem{subsemi}
J.~East, A.~{ Egri-Nagy}, and {J. D. Mitchell}, \emph{\textsc{{S}ub{S}emi} --
  GAP package for enumerating subsemigroups, v0.86}, 2025,
  \href{https://gap-packages.github.io/subsemi/}{\url{https://gap-packages.github.io/subsemi/}}.

\bibitem{T3enum1970}
\BIBentryALTinterwordspacing
C.~Wilde and S.~Raney, ``Computation of the transformation semigroups on three
  letters,'' \emph{Journal of the Australian Mathematical Society}, vol.~14,
  pp. 335--335, 11 1972. [Online]. Available:
  \url{http://journals.cambridge.org/article_S1446788700010806}
\BIBentrySTDinterwordspacing

\bibitem{T3enum1991}
\BIBentryALTinterwordspacing
G.~Bijev and K.~Todorov, ``\BIBforeignlanguage{English}{On the representation
  of abstract semigroups by transformation semigroups: Computer
  investigations},'' \emph{\BIBforeignlanguage{English}{Semigroup Forum}},
  vol.~43, no.~1, pp. 253--256, 1991. [Online]. Available:
  \url{http://dx.doi.org/10.1007/BF02574268}
\BIBentrySTDinterwordspacing

\bibitem{T4enum}
J.~East, A.~Egri-Nagy, and { J. D. Mitchell}, ``Enumerating transformation
  semigroups,'' \emph{Semigroup Forum}, vol.~95, no.~1, pp. 109--125, 2017.

\bibitem{diagsgps}
J.~East, A.~Egri-Nagy, A.~R. Francis, and J.~D. Mitchell, ``Finite diagram
  semigroups: Extending the computational horizon,'' 2015,
  \href{http://arxiv.org/abs/1502.07150}{arXiv:1502.07150 [math.GR]}.

\bibitem{SnIS}
A.~Egri-Nagy and V.~Gebhardt, ``Computational enumeration of independent
  generating sets of finite symmetric groups,'' 2016,
  \href{http://arxiv.org/abs/1602.03957}{arXiv:1602.03957 [math.GR]}.

\bibitem{GRAPEpaper}
L.~Soicher, ``Computing with graphs and groups,'' in \emph{Topics in Algebraic
  Graph Theory}.\hskip 1em plus 0.5em minus 0.4em\relax Cambridge University
  Press, 2004.

\bibitem{bliss}
T.~Junttila and P.~Kaski, ``\BIBforeignlanguage{English}{Engineering an
  efficient canonical labeling tool for large and sparse graphs},'' in
  \emph{\BIBforeignlanguage{English}{Proceedings of the ninth Workshop on
  Algorithm Engineering and Experiments and the fourth Workshop on Analytic
  Algorithmics and Combinatorics}}, 2007, pp. 135--149, society for Industrial
  and Applied Mathematics cop.; 978-0-898716-28-3; Workshop on Algorithm
  Engineering and Experiments.

\bibitem{nauty}
B.~D. McKay and A.~Piperno, ``Practical graph isomorphism, {II},''
  \emph{Journal of Symbolic Computation}, vol.~60, pp. 94--112, 2014.

\bibitem{Digraphs}
\BIBentryALTinterwordspacing
J.~D. Beule, J.~Jonu{\v s}as, J.~D. Mitchell, M.~Torpey, M.~Tsalakou, and W.~A.
  Wilson, ``Digraphs - {GAP} package, version 1.13.1,'' Sep 2025. [Online].
  Available: \url{https://digraphs.github.io/Digraphs}
\BIBentrySTDinterwordspacing

\bibitem{ComputingAutGrpOfSgps2010}
J.~Ara\'ujo, P.~B{\"u}nau, J.~Mitchell, and M.~Neunh{\"o}ffer, ``Computing
  automorphisms of semigroups,'' \emph{Journal of Symbolic Computation},
  vol.~45, no.~3, pp. 373 -- 392, 2010.

\bibitem{smallsemi}
A.~Distler and J.~Mitchell, ``{Smallsemi}, {A library of small semigroups},
  {V}ersion 0.7.2,''
  \href{https://gap-packages.github.io/smallsemi/}{\texttt{https://gap-packages.github.io/smallsemi}},
  Feb 2025, gAP package.

\bibitem{kigen}
A.~Egri-Nagy, ``{\textsc{kigen} Computational Semigroup Theory Software System
  written in {C}lojure},''
  \href{https://codeberg.org/egri-nagy/kigen}{\texttt{codeberg.org/egri-nagy/kigen}},
  2025.

\bibitem{2020Clojure}
\BIBentryALTinterwordspacing
R.~Hickey, ``A history of {C}lojure,'' \emph{Proc. ACM Program. Lang.}, vol.~4,
  no. HOPL, Jun. 2020. [Online]. Available:
  \url{https://doi.org/10.1145/3386321}
\BIBentrySTDinterwordspacing

\bibitem{Howie95}
J.~M. Howie, \emph{{Fundamentals of Semigroup Theory}}, ser. {London
  Mathematical Society Monographs New Series}.\hskip 1em plus 0.5em minus
  0.4em\relax Oxford University Press, 1995, vol.~12.

\bibitem{2025transformationrepresentationsdiagrammonoids}
\BIBentryALTinterwordspacing
R.~Cirpons, J.~East, and J.~D. Mitchell, ``Transformation representations of
  diagram monoids,'' 2025. [Online]. Available:
  \url{https://arxiv.org/abs/2411.14693}
\BIBentrySTDinterwordspacing

\bibitem{2014JEPartitionMonoidEmbeddings}
J.~East, ``\BIBforeignlanguage{English}{Partition monoids and embeddings in
  2-generator regular *-semigroups},''
  \emph{\BIBforeignlanguage{English}{Periodica Mathematica Hungarica}},
  vol.~69, no.~2, pp. 211--221, 2014.

\bibitem{WelshPowell67}
\BIBentryALTinterwordspacing
D.~J.~A. Welsh and M.~B. Powell, ``An upper bound for the chromatic number of a
  graph and its application to timetabling problems,'' \emph{The Computer
  Journal}, vol.~10, no.~1, pp. 85--86, 01 1967. [Online]. Available:
  \url{https://doi.org/10.1093/comjnl/10.1.85}
\BIBentrySTDinterwordspacing

\end{thebibliography}

\end{document}